\documentclass[12pt]{article}

\usepackage{amsmath}
\usepackage{amsfonts}
\usepackage{graphicx}
\usepackage{authblk}
\usepackage{multirow}

\newtheorem{lemma}{LEMMA}[section]
\newtheorem{theorem}{THEOREM}[section]
\newtheorem{corollary}{COROLLARY}[section]

\newtheorem{algorithm}{ALGORITHM}[section]
\newtheorem{example}{EXAMPLE}[section]

\newcommand{\ds}{\displaystyle}
\newcommand{\be}{\begin{equation}}
\newcommand{\ee}{\end{equation}}
\newcommand{\bp}{\underline{\bf Proof}:\ }
\newcommand{\ep}{{\hfill $\Box$}\\ }
\newcommand{\nt}{\noindent}

\title{On the Convergence of the Accelerated Riccati Iteration Method}

\author[1]{Prasanthan Rajasingam\footnote{E-mail: {\tt prasanthanr@jfn.ac.lk}}}
\author[2]{Jianhong Xu\footnote{Corresponding author. E-mail: {\tt jhxu@siu.edu}, 
Phone: +1 (618) 453--6510, Fax: +1 (618) 453--5300. The research of this author is partially supported by NSF Grant DMS-1419028.}}

\affil[1]{Department of Mathematics \& Statistics, University of Jaffna, Jaffna, Sri Lanka}
\affil[2]{Department of Mathematics, Southern Illinois University Carbondale, Carbondale, Illinois 62901, U.S.A.}
\begin{document}

\maketitle

\begin{abstract}
In this paper, we establish results fully addressing two open problems proposed recently by I. Ivanov, see Nonlinear Analysis 69 (2008) 4012--4024, with respect to the convergence of the accelerated Riccati iteration method for solving the continuous coupled algebraic Riccati equation, or CCARE for short. These results confirm several desirable features of that method, including the monotonicity and boundedness of the sequences it produces, its capability of determining whether the CCARE has a solution, the extremal solutions it computes under certain circumstances, and its faster convergence than the regular Riccati iteration method.
\end{abstract}

\nt {\bf Keywords}: continuous coupled algebraic Riccati equation, Markovian jump linear system, accelerated iteration method, convergence, rate of convergence, monotonicity, positive semidefinite solution, extremal solution

\nt {\bf AMS Subject Classification}: 15A24, 15B48, 34H05, 65B99, 65F30, 93B40

%\begin{keywords}
%\end{keywords}

%\begin{AMS}
%\end{AMS}

%\pagestyle{myheadings}
%\thispagestyle{plain}
%\markboth{J. Xu}{}

\section{Introduction}
\label{intro}

In this paper, all matrices are real and square. The size of a matrix may not be specified if it is clear from the context. For the sake of brevity, a positive semidefinite matrix $X$ is denoted by $X \succeq 0$. The term positive semidefinite, by convention, refers here only to the symmetric case, namely $X^T=X$. For symmetric matrices $X$ and $Y$, $X \succeq Y$ means $X-Y \succeq 0$. Similarly, $X \preceq Y$ means $Y-X \succeq 0$. In addition, $\langle N \rangle$ stands for $\{1, 2, \ldots, N\}$.
\\

The main problems we shall address in this paper concern the so-called continuous coupled algebraic Riccati equation, abbreviated as CCARE from now on. Specifically, let $A_i, S_i, Q_i \in \mathbb R^{n \times n}$, where $i \in \langle N \rangle$, and suppose that $S_i \succeq 0$ and $Q_i \succeq 0$ for all $i$, then the CCARE can be expressed in the form \cite{Guo13, Iva08}
\be
\label{ccare}
A_i^TX_i+X_iA_i-X_iS_iX_i+\sum_{j \in \langle N \rangle \backslash \{i\}}\delta_{i,j}X_j+Q_i=0, ~i \in \langle N \rangle,
\ee
where $\delta_{i,j} \ge 0$ for any $i \ne j$ and, moreover, $\ds \sum_{j \in \langle N \rangle \backslash \{i\}}\delta_{i,j}>0$ for each $i$. When there is no ambiguity, we shall denote by $X_i$, with $i \in \langle N \rangle$, a solution to the CCARE and call each $X_i$ the $i$th component of the solution.
\\

In particular, when $N=1$, the CCARE reduces to the classical continuous algebraic Riccati equation, or CARE for short in the sequel, which can be written by removing the subscript $i$ as 
\be
\label{care}
A^TX+XA-XSX+Q=0,
\ee
where $S \succeq 0$ and $Q \succeq 0$. Throughout this paper, we shall always assume by default that $N \ge 2$ in (\ref{ccare}). The CARE in (\ref{care}), however, plays a critical role in dealing with the main problems here.
\\

The CCARE in (\ref{ccare}) arises originally from an optimal control problem on Markovian jump linear systems. For background material, see, for example, \cite{CosFraTod13,Mar90}. Due to its connection to the solution of the optimal control problem, the numerical computation of positive semidefinite solutions to the CCARE has drawn much attention in the literature, see \cite{AboFreJan94, CosVal04, DamHin01, GajBor95, Guo13, Iva07, Iva08, ValGerCos99} along with the references therein. Among these, the following two numerical methods are relevant here: one is the Riccati iteration method, while the other is the accelerated (or modified) Riccati iteration method.
\\

We recall in passing the concepts of stabilizability and detectability. Let $A, S, Q \in \mathbb R^{n \times n}$. Then, $(A, S)$ is called stabilizable if there exists matrix $K$ such that $A-SK$ is stable, whereas $(A, Q)$ is called detectable if $(A^T, Q^T)$ is stabilizable. As a well-known result, such conditions guarantee the existence and uniqueness of a positive semidefinite solution to the CARE. This result will be stated formally in the next section.
\\

The Riccati iteration method and its convergence are investigated in \cite{CosVal04}. This method can be formulated --- see also \cite{Iva08} --- as:

\begin{algorithm}
\label{ricca}
For each $i \in \langle N \rangle$, choose the initial $X_i^{(0)} \succeq 0$ and set $\rho_i \ge 0$ such that $(A_i-\rho_i I, S_i)$ is stabilizable and $(A_i-\rho_i I, Q_i)$ is detectable. Next, for $k=0, 1, 2, \ldots$, we iterate according to 
\be
\label{algo1}
\begin{array}{l}
(A_i-\rho_i I)^TX_i^{(k+1)}+X_i^{(k+1)}(A_i-\rho_i I)-X_i^{(k+1)}S_iX_i^{(k+1)}\\
\\
\ds +\sum_{j \in \langle N \rangle \backslash \{i\}}\delta_{i,j}X_j^{(k)}+Q_i+2\rho_i X_i^{(k)}=0, ~i \in \langle N \rangle.
\end{array}
\ee
\end{algorithm}

At each iteration, the above algorithm solves $N$ CARE's, either in serial or in parallel if all $X_i^{(k)}$'s are available, which may be implemented easily in practice with Matlab's {\tt care}.  As mentioned in \cite{CosVal04}, however, the main advantage of Algorithm \ref{ricca} is that the stabilizability and detectability conditions in this algorithm, i.e. in (\ref{algo1}), can always be satisfied by choosing appropriate values of $\rho_i$'s, and thus (\ref{algo1}) computes unique sequences of positive semidefinite matrices $\{X_i^{(k)}\}$, $i \in \langle N \rangle$, even when the CCARE in (\ref{care}) has no solution. Moreover, for each $i$, $\{X_i^{(k)}\}$ converges if and only if (\ref{care}) has a solution, and it does so in a monotonically increasing fashion toward the minimal solution of (\ref{care}), provided that $X_i^{(0)}=0$ for all $i$; see \cite{CosVal04} for more detail. Note that the latter feature here is especially attractive, since it means that the algorithm can also determine whether (\ref{care}) has a solution or not.
\\

The accelerated Riccati iteration method appears in \cite[(20)]{Iva08} as an effort to improve upon Algorithm \ref{ricca} via making use of updated $X_i^{(k+1)}$'s in (\ref{algo1}) as soon as they become available. Intuitively, such a modification should speed up the convergence of Algorithm \ref{ricca}. Specifically, this accelerated algorithm can be summarized as: 

\begin{algorithm}
\label{accel}
For each $i \in \langle N \rangle$, choose the initial $X_i^{(0)} \succeq 0$ and set $\rho_i \ge 0$ such that $(A_i-\rho_i I, S_i)$ is stabilizable and $(A_i-\rho_i I, Q_i)$ is detectable. Next, for $k=0, 1, 2, \ldots$, we iterate according to 
\be
\label{algo2}
\begin{array}{l}
(A_i-\rho_i I)^TX_i^{(k+1)}+X_i^{(k+1)}(A_i-\rho_i I)-X_i^{(k+1)}S_iX_i^{(k+1)}\\
\\
\ds +\sum_{j=1}^{i-1}\delta_{i,j}X_j^{(k+1)}+\sum_{j=i+1}^N\delta_{i,j}X_j^{(k)}+Q_i+2\rho_i X_i^{(k)}=0, ~i=1, 2, \ldots, N.
\end{array}
\ee
\end{algorithm}

Similar to the preceding one, at each iteration, the above accelerated algorithm solves $N$ CARE's, but clearly only in a serial fashion --- a potential trade-off between intrinsic parallelism and rate of convergence. Other shared features between the two algorithms are also expected here, such as the ease of implementation with available software and the existence and uniqueness of the sequences $\{X_i^{(k)}\}$, $i \in \langle N \rangle$, out of (\ref{algo2}), consisting entirely of positive semidefinite matrices. Nevertheless, Algorithm \ref{accel} poses a number of interesting and crucial problems too. Despite some favorable numerical evidence in \cite{Iva08}, the following questions remain yet to be explored \cite[p. 4021]{Iva08}:
\begin{itemize}
\item[]{\underline{Question 1}: What conditions are needed for (\ref{algo2}) to compute monotone, convergent sequences $\{X_i^{(k)}\}$, $i \in \langle N \rangle$?}
\item[]{\underline{Question 2}: Do such sequences converge faster in comparison to their counterparts from Algorithm \ref{ricca}?}
\end{itemize} 

The goals of this paper are to resolve these open problems that are vital to Algorithm \ref{accel}. 

\section{Convergence of Accelerated Riccati Iteration Method}
\label{conv}

Let us start with several necessary preparatory results on the solution of the CARE given by (\ref{care}).
\\

The first result here gives the necessary and sufficient conditions for the existence and uniqueness of a positive semidefinite solution to the CARE in terms of stabilizability and detectability.
\begin{lemma}{\rm (\cite[Theorem 2.21]{BinLanMei12})}
\label{lem1}
The CARE in (\ref{care}) has a unique positive semidefinite solution $X$ such that $A-SX$ is stable, namely $X$ is also stabilizing, if and only if $(A, S)$ is stabilizable and $(A, Q)$ is detectable.
\end{lemma}

The second result establishes an ordering for the solutions to (\ref{care}) under a varying term $Q$. For convenience of application, we reformulate it based on its original form in \cite{Wil71}.
\begin{lemma}{\rm (\cite[Lemma 3]{Wil71}, also \cite[Proposition 1]{CosVal04})}
\label{lem2}
Suppose that $S \succeq 0$ and $Q$ is symmetric. Let $X_1 \succeq 0$ be a solution of 
$$A^TX+XA-XSX+Q \preceq 0$$
such that $A-SX_1$ is stable and let $X_2 \succeq 0$ be a solution of 
$$A^TX+XA-XSX+Q \succeq 0.$$ Then, $X_1 \succeq X_2$.
\end{lemma}

Finally, we cite below a result concerning detectability. Its original proof in \cite{CosVal04} employs a rank argument, but it can also be shown alternatively using a well-known characterization of detectability.

\begin{lemma}{\rm (\cite[Proposition 2]{CosVal04})}
\label{lem3}
Suppose that $Q \succeq 0$ and $\Delta Q \succeq 0$. Then, $(A, Q+\Delta Q)$ is detectable whenever so is $(A, Q)$.
\end{lemma}
\bp
By the Popov-Belevitch-Hautus tests, see \cite[Theorem 8.5]{WilLaw07}, $(A, Q)$ is detectable if and only if there exists no (right) eigenvector $u$ of $A$ associated with eigenvalue $\lambda$ with ${\rm Re} \lambda \ge 0$ such that $Qu=0$.
\\

Let $(A, Q)$ be detectable. Suppose now to the contrary that $(A, Q+\Delta Q)$ is not detectable. We denote by $(\lambda, u)$, with ${\rm Re}\lambda \ge 0$, an eigenpair of $A$ such that $(Q+\Delta Q)u=0$. This leads to $u^\ast (Q+\Delta Q)u=u^\ast Qu+u^\ast \Delta Qu=0$. In particular, we have $u^\ast Qu=0$ and, consequently, $Qu=0$, which is a contradiction to the detectability of $(A, Q)$.
\ep

To facilitate the statement of our results, following \cite{Iva08}, we define that for each $i \in \langle N \rangle$, 
\be
\label{res}
{\cal R}_i(X_1, X_2, \ldots, X_N)=A_i^TX_i+X_iA_i-X_iS_iX_i+\sum_{j \in \langle N \rangle \backslash \{i\}}\delta_{i,j}X_j+Q_i.
\ee
Accordingly, the CCARE in (\ref{ccare}) can also be written as $${\cal R}_i(X_1, X_2, \ldots, X_N)=0, ~i \in \langle N \rangle.$$

We are now in a position to develop a number of results concerning the first question raised in \cite{Iva08}, i.e. sufficient conditions so as to guarantee that the accelerated Riccati iteration method in Algorithm \ref{accel} computes unique monotonically increasing, bounded sequences of positive semidefinite matrices $\{X_i^{(k)}\}$, $i \in \langle N \rangle$.

\begin{theorem}
\label{thm1}
Let $\hat X_i \succeq 0$, $i \in \langle N \rangle$, be such that for each $i$, $${\cal R}_i(\hat X_1, \hat X_2, \ldots, \hat X_N) \preceq 0.$$ In addition, suppose that the initial positive semidefinite $X_i^{(0)}$'s in Algorithm \ref{accel} are such that $${\cal R}_i(X_1^{(0)}, X_2^{(0)}, \ldots, X_N^{(0)}) \succeq 0$$ and $X_i^{(0)} \preceq \hat X_i$ for all $i \in \langle N \rangle$. Moreover, for each $i$, let $\rho_i \ge 0$ be such that $(A_i-\rho_i I, S_i)$ is stabilizable, $(A_i-\rho_i I, Q_i)$ is detectable, and $A_i-\rho_i I-S_i\hat X_i$ is stable. Then,
\begin{itemize}
\item[(i)]{Algorithm \ref{accel} computes unique sequences of positive semidefinite matrices $\{X_i^{(k+1)}\}$, where $k=0, 1, 2, \ldots$ and $i \in \langle N \rangle$. Besides, for each $i$, $A_i-\rho_i I -S_iX_i^{(k+1)}$, where $k=0, 1, 2, \ldots$, are all stable.}
\item[(ii)]{For each $i$, $X_i^{(k+1)} \succeq X_i^{(k)}$ for all $k=0, 1, 2, \ldots$; that is, each $\{X_i^{(k)}\}$ is monotonically increasing.}
\item[(iii)]{For each $i$, ${\cal R}_i(X_1^{(k)}, X_2^{(k)}, \ldots, X_N^{(k)}) \succeq 0$ for all $k=0, 1, 2, \ldots$.}
\item[(iv)]{For each $i$, $X_i^{(k)} \preceq \hat X_i$ for all $k=0, 1, 2, \ldots$; that is, each $\{X_i^{(k)}\}$ is also bounded above.}
\end{itemize}
\end{theorem}
\bp
We proceed by way of induction on $k$ and, for each $k$, induction on $i$ as well.
\\

\underline{Case $k=0$}: In this case, (iii) and (iv) are trivially true by assumption.
\\

Let $i=1$. From (\ref{algo2}), we have 
\be
\label{eqn1}
(A_1-\rho_1 I)^TX_1^{(1)}+X_1^{(1)}(A_1-\rho_1 I)-X_1^{(1)}S_1X_1^{(1)}+Q_1+\Delta Q_1=0,
\ee
where $\ds \Delta Q_1=\sum_{j=2}^N\delta_{1,j}X_j^{(0)}+2\rho_1 X_1^{(0)} \succeq 0$. Since $(A_1-\rho_1 I, S_1)$ is stabilizable and, following Lemma \ref{lem3}, $(A_1-\rho_1 I, Q_1+\Delta Q_1)$ is detectable, we know by Lemma \ref{lem1} that (\ref{eqn1}) has a unique solution $X_1^{(1)} \succeq 0$ such that $A_1-\rho_1 I -S_1X_1^{(1)}$ is stable, and hence (i) holds at $k=0$ and $i=1$. In addition, using ${\cal R}_1(X_1^{(0)}, X_2^{(0)}, \ldots, X_N^{(0)}) \succeq 0$, we have 
\be
\label{eqn2}
(A_1-\rho_1 I)^TX_1^{(0)}+X_1^{(0)}(A_1-\rho_1 I)-X_1^{(0)}S_1X_1^{(0)}+Q_1+\Delta Q_1 \succeq 0,
\ee
where $\Delta Q_1$ is given as below (\ref{eqn1}). It follows from (\ref{eqn1}), (\ref{eqn2}), the stability of $A_1-\rho_1 I -S_1X_1^{(1)}$, and Lemma \ref{lem2} that $X_1^{(1)} \succeq X_1^{(0)}$, i.e. (ii) holds as well at $k=0$ and $i=1$.
\\

Suppose next that for some $2 \le r \le N$, (i) and (ii) are justified at $k=0$ for all $i=1, 2, \ldots, r-1$. On letting $i=r$ in (\ref{algo2}), we obtain 
\be
\label{eqn3}
(A_r-\rho_r I)^TX_r^{(1)}+X_r^{(1)}(A_r-\rho_r I)-X_r^{(1)}S_rX_r^{(1)}+Q_r+\Delta Q_r=0,
\ee
where $\ds \Delta Q_r=\sum_{j=1}^{r-1}\delta_{r,j}X_j^{(1)}+\sum_{j=r+1}^N\delta_{r,j}X_j^{(0)}+2\rho_rX_r^{(0)} \succeq 0$. Since $(A_r-\rho_r I, S_r)$ is stabilizable while, from Lemma \ref{lem3}, $(A_r-\rho_r I, Q_r+\Delta Q_r)$ is detectable, (\ref{eqn3}) has a unique solution $X_r^{(1)} \succeq 0$ such that $A_r-\rho_r I -S_rX_r^{(1)}$ is stable according to Lemma \ref{lem1}, and hence (i) holds at $k=0$. Finally, observe that ${\cal R}_r(X_1^{(0)}, X_2^{(0)}, \ldots, X_N^{(0)}) \succeq 0$ and $X_i^{(1)} \succeq X_i^{(0)}$, $i=1, 2, \ldots, r-1$, yield 
\be
\label{eqn4}
(A_r-\rho_r I)^TX_r^{(0)}+X_r^{(0)}(A_r-\rho_r I)-X_r^{(0)}S_rX_r^{(0)}+Q_r+\Delta Q_r \succeq 0,
\ee
where $\Delta Q_r$ is given under (\ref{eqn3}). Due to (\ref{eqn3}), (\ref{eqn4}), the stability of $A_r-\rho_r I -S_rX_r^{(1)}$, and Lemma \ref{lem2}, we see that $X_r^{(1)} \succeq X_r^{(0)}$, i.e. (ii) holds too at $k=0$.
\\

This concludes the proof of (i) through (iv) for the case $k=0$.
\\

\underline{Case $k > 0$}: Suppose now that (i) through (iv) are true for some $k \ge 0$. We show here that they remain true at $k+1$.
\\

First, by (\ref{algo2}), and with (ii) and (iii) being true at $k$, it is clear that $${\cal R}_i(X_1^{(k+1)}, X_2^{(k+1)}, \ldots, X_N^{(k+1)}) \succeq 0, ~i \in \langle N \rangle,$$ i.e. (iii) holds at $k+1$. 
\\

Next, for (i) and (ii), we start with $i=1$. Using (\ref{algo2}), we have 
\be
\label{eqn5}
(A_1-\rho_1 I)^TX_1^{(k+2)}+X_1^{(k+2)}(A_1-\rho_1 I)-X_1^{(k+2)}S_1X_1^{(k+2)}+Q_1+\Delta \tilde Q_1=0,
\ee
where $\ds \Delta \tilde Q_1=\sum_{j=2}^N\delta_{1,j}X_j^{(k+1)}+2\rho_1X_1^{(k+1)} \succeq 0$. Since $(A_1-\rho_1 I, S_1)$ is stabilizable and, via Lemma \ref{lem3}, $(A_1-\rho_1 I, Q_1+\Delta \tilde Q_1)$ is detectable, in view of Lemma \ref{lem1}, (\ref{eqn5}) has a unique solution $X_1^{(k+2)} \succeq 0$ with $A_1-\rho_1 I -S_1X_1^{(k+2)}$ being stable and, consequently, (i) is true at $k+1$ and $i=1$. In addition, we find from ${\cal R}_1(X_1^{(k+1)}, X_2^{(k+1)}, \ldots, X_N^{(k+1)}) \succeq 0$ that 
\be
\label{eqn6}
(A_1-\rho_1 I)^TX_1^{(k+1)}+X_1^{(k+1)}(A_1-\rho_1 I)-X_1^{(k+1)}S_1X_1^{(k+1)}+Q_1+\Delta \tilde Q_1 \succeq 0,
\ee
where $\Delta \tilde Q_1$ is given following (\ref{eqn5}). Because of the stability of $A_1-\rho_1 I - S_1X_1^{(k+2)}$ and Lemma \ref{lem2}, (\ref{eqn5}), and (\ref{eqn6}) imply $X_1^{(k+2)} \succeq X_1^{(k+1)}$, showing that (ii) also holds at $k+1$ and $i=1$.
\\

Suppose now that for some $2 \le r \le N$, both (i) and (ii) hold true for $i=1, 2, \ldots, r-1$ at $k+1$. According to (\ref{algo2}), we have 
\be
\label{eqn7}
(A_r-\rho_r I)^TX_r^{(k+2)}+X_r^{(k+2)}(A_r-\rho_r I)-X_r^{(k+2)}S_rX_r^{(k+2)}+Q_r+\Delta \tilde Q_r=0,
\ee
where $\ds \Delta \tilde Q_r=\sum_{j=1}^{r-1}\delta_{r,j}X_j^{(k+2)}+\sum_{j=r+1}^N\delta_{r,j}X_j^{(k+1)}+2\rho_rX_r^{(k+1)} \succeq 0$. Observe that, from Lemma \ref{lem3}, $(A_r-\rho_r I, Q_r+\Delta \tilde Q_r)$ is detectable. Besides, $(A_r-\rho_r I, S_r)$ is stabilizable. Hence, by Lemma \ref{lem1}, (\ref{eqn7}) has a unique solution $X_r^{(k+2)} \succeq 0$ such that $A_r-\rho_r I -S_rX_r^{(k+2)}$ is stable, implying that (i) is true at $k+1$. Finally, combining ${\cal R}_r(X_1^{(k+1)}, X_2^{(k+1)}, \ldots, X_N^{(k+1)}) \succeq 0$ and $X_i^{(k+2)} \succeq X_i^{(k+1)}$, $i=1, 2, \ldots, r-1$, we arrive at 
\be
\label{eqn8}
(A_r-\rho_r I)^TX_r^{(k+1)}+X_r^{(k+1)}(A_r-\rho_r I)-X_r^{(k+1)}S_rX_r^{(k+1)}+Q_r+\Delta \tilde Q_r \succeq 0,
\ee
where $\Delta \tilde Q_r$ is given under (\ref{eqn7}). By Lemma \ref{lem2}, (\ref{eqn7}), (\ref{eqn8}), and the stability of $A_r-\rho_r I -S_rX_r^{(k+2)}$ yield $X_r^{(k+2)} \succeq X_r^{(k+1)}$, i.e. (ii) holds at $k+1$.
\\

It remains to show that (iv) is true at $k+1$, i.e. $X_i^{(k+1)} \preceq \hat X_i$, $i \in \langle N \rangle$. Again, we start with $i=1$. On one hand, because of ${\cal R}_1(\hat X_1, \hat X_2, \ldots, \hat X_N) \preceq 0$, we have 
\be
\label{eqn9}
(A_1-\rho_1 I)^T\hat X_1+\hat X_1(A_1-\rho_1 I)-\hat X_1S_1\hat X_1+Q_1+\Delta \bar Q_1 \preceq 0,
\ee
where $\ds \Delta \bar Q_1 = \sum_{j=2}^N\delta_{1,j}\hat X_j+2\rho_1\hat X_1 \succeq 0$. On the other hand, seeing (\ref{algo2}) along with $X_i^{(k)} \preceq \hat X_i$, $i \in \langle N \rangle$, we obtain 
\be
\label{eqn10}
(A_1-\rho_1 I)^TX_1^{(k+1)}+X_1^{(k+1)}(A_1-\rho_1 I)-X_1^{(k+1)}S_1X_1^{(k+1)}+Q_1+\Delta \bar Q_1 \succeq 0,
\ee
where $\Delta \bar Q_1$ is given as below (\ref{eqn9}). Since $A_1-\rho_1 I - S_1\hat X_1$ is stable, accordingly to Lemma \ref{lem2}, we get from (\ref{eqn9}) and (\ref{eqn10}) that $X_1^{(k+1)} \preceq \hat X_1$.
\\

Suppose next that for some $2 \le r \le N$, $X_i^{(k+1)} \preceq \hat X_i$, $i=1, 2, \ldots, r-1$. By ${\cal R}_r(\hat X_1, \hat X_2, \ldots, \hat X_N) \preceq 0$, we get 
\be
\label{eqn11}
(A_r-\rho_r I)^T\hat X_r+\hat X_r(A_r-\rho_r I)-\hat X_rS_r\hat X_r+Q_r+\Delta \bar Q_r \preceq 0,
\ee
where $\ds \Delta \bar Q_r=\sum_{j=1}^{r-1}\delta_{r,j}\hat X_j + \sum_{j=r+1}^N\delta_{r,j}\hat X_j + 2\rho_r\hat X_r \succeq 0$. In the meantime, we use (\ref{algo2}) together with $X_i^{(k+1)} \preceq \hat X_i$, where $i=1, 2, \ldots, r-1$, and $X_i^{(k)} \preceq \hat X_i$, where $i=r+1, r+2, \ldots, N$, to derive 
\be
\label{eqn12}
(A_r-\rho_r I)^TX_r^{(k+1)}+X_r^{(k+1)}(A_r-\rho_r I)-X_r^{(k+1)}S_rX_r^{(k+1)}+Q_r+\Delta \bar Q_r \succeq 0,
\ee
where $\Delta \bar Q_r$ is given next to (\ref{eqn11}). Finally, the stability of $A_r-\rho_r I - S_r\hat X_r$, (\ref{eqn11}), (\ref{eqn12}), and Lemma \ref{lem2} lead to $X_r^{(k+1)} \preceq \hat X_r$. This shows that (iv) holds too at $k+1$.
\\

The proof is now complete in its entirety.
\ep

An immediate consequence of Theorem \ref{thm1} goes as follows.

\begin{corollary}
\label{cor1}
Under the same conditions as Theorem \ref{thm1}, Algorithm \ref{accel} computes unique sequences of positive semidefinite matrices $\{X_i^{(k)}\}$, with $i \in \langle N \rangle$, that converge to a positive semidefinite solution $X_i$, $i \in \langle N \rangle$, of the CCARE in (\ref{ccare}), i.e. $\ds \lim_{k \rightarrow \infty}X_i^{(k)}=X_i$ for each $i$.
\end{corollary}
\bp
For each $i$, the convergence of $\{X_i^{(k)}\}$ is obvious --- see, for example, \cite{IvaHasMin01} and \cite[Corollary 4.1]{XuXia13} --- and it does so toward some positive semidefinite $X_i$. Next, by pushing $k \rightarrow \infty$ in (\ref{algo2}), we see that $X_i$, $i \in \langle N \rangle$, is indeed a solution to (\ref{ccare}).
\ep

Corollary \ref{cor1} shows that, similar to the pure Riccati iteration method in Algorithm \ref{ricca}, the accelerated version here in Algorithm \ref{accel} can also determine whether the CCARE has a solution or not. To be specific, Algorithm \ref{accel} yields a positive semidefinite solution to the CCARE whenever it converges. 
\\

Furthermore, if $\hat X_i$'s in Theorem \ref{thm1} happen to be a positive semidefinite solution to the CCARE in (\ref{ccare}), then Algorithm \ref{accel} actually finds the minimal positive semidefinite solution to (\ref{ccare}) as the next result demonstrates.

\begin{corollary}
\label{cor2}
Let $X_i \succeq 0$, $i \in \langle N \rangle$, be a solution to (\ref{ccare}). Suppose that the initial positive semidefinite $X_i^{(0)}$'s in Algorithm \ref{accel} are such that $${\cal R}_i(X_1^{(0)}, X_2^{(0)}, \ldots, X_N^{(0)}) \succeq 0$$ and $X_i^{(0)} \preceq X_i$ for all $i \in \langle N \rangle$. Moreover, for each $i$, let $\rho_i \ge 0$ be such that $(A_i-\rho_i I, S_i)$ is stabilizable, $(A_i-\rho_i I, Q_i)$ is detectable, and $A_i-\rho_i I - S_iX_i$ is stable. Then, Algorithm \ref{accel} produces unique sequences of positive semidefinite matrices $\{X_i^{(k)}\}$, $i \in \langle N \rangle$, such that for each $i$, $\{X_i^{(k)}\}$ is monotonically increasing, bounded above by $X_i$, and converges to $X_i^-$, the $i$th component of the minimal positive semidefinite solution to the CCARE in (\ref{ccare}).
\end{corollary}
\bp
It is clear that Corollary \ref{cor2} assumes the same conditions as Theorem \ref{thm1}, except for $\hat X_i$'s being replaced with $X_i$'s. By Corollary \ref{cor1}, we know that Algorithm \ref{accel} computes a positive semidefinite solution $X_i^-$, $i \in \langle N \rangle$, to (\ref{ccare}). Besides, note that due to Theorem \ref{thm1}, $X_i^- \preceq X_i$ for all $i$ whenever $X_i$, $i \in \langle N \rangle$, is a solution to (\ref{ccare}), thus $X_i^-$, $i \in \langle N \rangle$, is the minimal positive semidefinite solution to (\ref{ccare}).   
\ep

Clearly, Corollaries \ref{cor1} and \ref{cor2} also verify that Algorithm \ref{accel} converges if and only if (\ref{ccare}) has a positive semidefinite solution, provided that the initial $X_i^{(0)}$'s are chosen as in Corollary \ref{cor2}. We point out that, in particular, those conditions on $X_i^{(0)}$'s are trivially satisfied when $X_i^{(0)}=0$, $i \in \langle N \rangle$. In other words, this desirable feature of the Riccati iteration method for allowing a full determination of the existence of a positive semidefinite solution --- see \cite{CosVal04} --- carries over to the accelerated Riccati iteration method here.
\\

In light of Theorem \ref{thm1} and Corollaries \ref{cor1} and \ref{cor2}, we can formulate the following three parallel results, whose proofs are very similar and, therefore, are omitted.

\begin{theorem}
\label{thm2}
Let $\hat X_i \succeq 0$, $i \in \langle N \rangle$, be such that for each $i$, $${\cal R}_i(\hat X_1, \hat X_2, \ldots, \hat X_N) \succeq 0.$$ In addition, suppose that the initial positive semidefinite $X_i^{(0)}$'s in Algorithm \ref{accel} are such that $${\cal R}_i(X_1^{(0)}, X_2^{(0)}, \ldots, X_N^{(0)}) \preceq 0$$ and $X_i^{(0)} \succeq \hat X_i$ for all $i \in \langle N \rangle$. Moreover, for each $i$, let $\rho_i \ge 0$ be such that $(A_i-\rho_i I, S_i)$ is stabilizable, $(A_i-\rho_i I, Q_i)$ is detectable, and $A_i-\rho_i I - S_iX_i^{(0)}$ is stable. Then,
\begin{itemize}
\item[(i)]{Algorithm \ref{accel} computes unique sequences of positive semidefinite matrices $\{X_i^{(k+1)}\}$, where $k=0, 1, 2, \ldots$ and $i \in \langle N \rangle$. Besides, for each $i$, $A_i-\rho_i I -S_iX_i^{(k+1)}$, where $k=0, 1, 2, \ldots$, are all stable.}
\item[(ii)]{For each $i$, $X_i^{(k+1)} \preceq X_i^{(k)}$ for all $k=0, 1, 2, \ldots$; that is, each $\{X_i^{(k)}\}$ is monotonically decreasing.}
\item[(iii)]{For each $i$, ${\cal R}_i(X_1^{(k)}, X_2^{(k)}, \ldots, X_N^{(k)}) \preceq 0$ for all $k=0, 1, 2, \ldots$.}
\item[(iv)]{For each $i$, $X_i^{(k)} \succeq \hat X_i$ for all $k=0, 1, 2, \ldots$; that is, each $\{X_i^{(k)}\}$ is also bounded below.}
\end{itemize}
\end{theorem}

\begin{corollary}
\label{cor3}
Under the same conditions as Theorem \ref{thm2}, Algorithm \ref{accel} computes unique sequences of positive semidefinite matrices $\{X_i^{(k)}\}$, with $i \in \langle N \rangle$, that converge to a positive semidefinite solution $X_i$, $i \in \langle N \rangle$, of the CCARE in (\ref{ccare}), i.e. $\ds \lim_{k \rightarrow \infty}X_i^{(k)}=X_i$ for each $i$.
\end{corollary}

\begin{corollary}
\label{cor4}
Let $X_i \succeq 0$, $i \in \langle N \rangle$, be a solution to (\ref{ccare}). Suppose that the initial positive semidefinite $X_i^{(0)}$'s in Algorithm \ref{accel} are such that $${\cal R}_i(X_1^{(0)}, X_2^{(0)}, \ldots, X_N^{(0)}) \preceq 0$$ and $X_i^{(0)} \succeq X_i$ for all $i \in \langle N \rangle$. Moreover, for each $i$, let $\rho_i \ge 0$ be such that $(A_i-\rho_i I, S_i)$ is stabilizable, $(A_i-\rho_i I, Q_i)$ is detectable, and $A_i-\rho_i I -S_iX_i^{(0)}$ is stable. Then, Algorithm \ref{accel} produces unique sequences of positive semidefinite matrices $\{X_i^{(k)}\}$, $i \in \langle N \rangle$, such that for each $i$, $\{X_i^{(k)}\}$ is monotonically decreasing, bounded below by $X_i$, and converges to $X_i^+$, the $i$th component of the maximal positive semidefinite solution to the CCARE in (\ref{ccare}).
\end{corollary}

Corollaries \ref{cor1} through \ref{cor4}, coupled with Theorems \ref{thm1} and \ref{thm2}, serve as a rather complete answer to the first open problem in \cite{Iva08}. Especially, these corollaries spell out not only the conditions for convergence in Algorithm \ref{accel} but also the particular extremal types of solution this algorithm converges to under certain circumstances.
\\

It is straightforward to see that, in fact, Algorithm \ref{ricca} shares all of the preceding results on Algorithm \ref{accel}. The proofs are very similar except that the inductive steps on $i$ are no longer needed. For the sake of concision, we only state such results without proof in forms parallel to Theorems \ref{thm1} and \ref{thm2}. In addition, for clarity, we denote the sequences from Algorithm \ref{ricca} by $\{Y_i^{(k)}\}$'s here.

\begin{theorem}
\label{thm3}
Let $\hat X_i \succeq 0$, $i \in \langle N \rangle$, be such that for each $i$, $${\cal R}_i(\hat X_1, \hat X_2, \ldots, \hat X_N) \preceq 0.$$ In addition, suppose that the initial positive semidefinite $Y_i^{(0)}$'s in Algorithm \ref{ricca} are such that $${\cal R}_i(Y_1^{(0)}, Y_2^{(0)}, \ldots, Y_N^{(0)}) \succeq 0$$ and $Y_i^{(0)} \preceq \hat X_i$ for all $i \in \langle N \rangle$. Moreover, for each $i$, let $\rho_i \ge 0$ be such that $(A_i-\rho_i I, S_i)$ is stabilizable, $(A_i-\rho_i I, Q_i)$ is detectable, and $A_i-\rho_i I-S_i\hat X_i$ is stable. Then,
\begin{itemize}
\item[(i)]{Algorithm \ref{ricca} computes unique sequences of positive semidefinite matrices $\{Y_i^{(k+1)}\}$, where $k=0, 1, 2, \ldots$ and $i \in \langle N \rangle$. Besides, for each $i$, $A_i-\rho_i I -S_iY_i^{(k+1)}$, where $k=0, 1, 2, \ldots$, are all stable.}
\item[(ii)]{For each $i$, $Y_i^{(k+1)} \succeq Y_i^{(k)}$ for all $k=0, 1, 2, \ldots$; that is, each $\{Y_i^{(k)}\}$ is monotonically increasing.}
\item[(iii)]{For each $i$, ${\cal R}_i(Y_1^{(k)}, Y_2^{(k)}, \ldots, Y_N^{(k)}) \succeq 0$ for all $k=0, 1, 2, \ldots$.}
\item[(iv)]{For each $i$, $Y_i^{(k)} \preceq \hat X_i$ for all $k=0, 1, 2, \ldots$; that is, each $\{Y_i^{(k)}\}$ is also bounded above.}
\end{itemize}
\end{theorem}

\begin{theorem}
\label{thm4}
Let $\hat X_i \succeq 0$, $i \in \langle N \rangle$, be such that for each $i$, $${\cal R}_i(\hat X_1, \hat X_2, \ldots, \hat X_N) \succeq 0.$$ In addition, suppose that the initial positive semidefinite $Y_i^{(0)}$'s in Algorithm \ref{ricca} are such that $${\cal R}_i(Y_1^{(0)}, Y_2^{(0)}, \ldots, Y_N^{(0)}) \preceq 0$$ and $Y_i^{(0)} \succeq \hat X_i$ for all $i \in \langle N \rangle$. Moreover, for each $i$, let $\rho_i \ge 0$ be such that $(A_i-\rho_i I, S_i)$ is stabilizable, $(A_i-\rho_i I, Q_i)$ is detectable, and $A_i-\rho_i I - S_iY_i^{(0)}$ is stable. Then,
\begin{itemize}
\item[(i)]{Algorithm \ref{ricca} computes unique sequences of positive semidefinite matrices $\{Y_i^{(k+1)}\}$, where $k=0, 1, 2, \ldots$ and $i \in \langle N \rangle$. Besides, for each $i$, $A_i-\rho_i I -S_iY_i^{(k+1)}$, where $k=0, 1, 2, \ldots$, are all stable.}
\item[(ii)]{For each $i$, $Y_i^{(k+1)} \preceq Y_i^{(k)}$ for all $k=0, 1, 2, \ldots$; that is, each $\{Y_i^{(k)}\}$ is monotonically decreasing.}
\item[(iii)]{For each $i$, ${\cal R}_i(Y_1^{(k)}, Y_2^{(k)}, \ldots, Y_N^{(k)}) \preceq 0$ for all $k=0, 1, 2, \ldots$.}
\item[(iv)]{For each $i$, $Y_i^{(k)} \succeq \hat X_i$ for all $k=0, 1, 2, \ldots$; that is, each $\{Y_i^{(k)}\}$ is also bounded below.}
\end{itemize}
\end{theorem}

Compared with the results in \cite{CosVal04}, Theorems \ref{thm3} and \ref{thm4} on Algorithm \ref{ricca} are broader because, firstly, they allow nonzero initial $Y_i^{(0)}$'s and, secondly, they provide respective sufficient conditions for the resulting convergent sequences $\{Y_i^{(k)}\}$'s to be either monotonically increasing or monotonically decreasing. Consequently, conclusions on extremal solutions Algorithm \ref{ricca} can compute follow easily from these theorems --- in a way similar to Corollaries \ref{cor2} and \ref{cor4}.
\\

We comment that in Theorem \ref{thm1}, Corollary \ref{cor2}, and Theorem \ref{thm3}, as in \cite{CosVal04, Iva08}, an easy choice of the initial $X_i^{(0)}$'s and $Y_i^{(0)}$'s is to set $X_i^{(0)}=Y_i^{(0)}=0$ for any $i \in \langle N \rangle$. With this choice, all the conditions on $X_i^{(0)}$'s and $Y_i^{(0)}$'s in those results are trivially satisfied. On the other hand, when applying Theorem \ref{thm2}, Corollary \ref{cor4}, and Theorem \ref{thm4}, we may choose the initial $X_i^{(0)}$'s and $Y_i^{(0)}$'s to be some existing upper solution bounds for the CCARE. For results relevant to such bounds, see, for example, \cite{CzoSwi01, DavShiWil08, Xu13, XuRaj16} and the references therein.
\\

One of the advantages shared by Algorithms \ref{ricca} and \ref{accel} is that the stabilizability and detectability requirements can always be met by appropriate values of $\rho_i$'s. In Theorem \ref{thm1}, Corollary \ref{cor2}, and Theorem \ref{thm3}, however, the choice of $\rho_i$'s is complicated by the stability requirement on $A_i-\rho_i I - S_iX_i$'s since, in practice, the solution $X_i$'s is not available {\it a priori}. Although this issue might be alleviated by resorting to sufficiently large $\rho_i$ values, we shall demonstrate later that, similar to Algorithm \ref{ricca}, unnecessarily large $\rho_i$ values are usually not advisable for Algorithm \ref{accel}.  
\\

Next, we move on to examining the other open problem in \cite{Iva08} regarding a comparison of the rate of convergence of the accelerated Riccati iteration method versus that of the Riccati iteration method. In this regard, we prove the following:  

\begin{theorem}
\label{thm5}
Under the same assumptions as in Theorems \ref{thm1} and \ref{thm3} with $X_i^{(0)}=Y_i^{(0)}$ for all $i \in \langle N \rangle$, on letting $\{X_i^{(k)}\}$, $i \in \langle N \rangle$, be the sequences computed with Algorithm \ref{accel} and $\{Y_i^{(k)}\}$, $i \in \langle N \rangle$, be the corresponding sequences computed with Algorithm \ref{ricca}, we have that for each $i$, $X_i^{(k)} \succeq Y_i^{(k)}$, where $k=0, 1, 2, \ldots$.
\end{theorem}
\bp
Again, we use induction on $k$ and, given $k$, induction on $i$. The case $k=0$ is trivial here.
\\

Suppose now that at some $k \ge 0$,
\be
\label{eqn13}
X_i^{(k)} \succeq Y_i^{(k)}, ~i \in \langle N \rangle.
\ee
Let us show that $X_i^{(k+1)} \succeq Y_i^{(k+1)}$, $i \in \langle N \rangle$.
\\

From (\ref{algo2}) and (\ref{eqn13}), we obtain 
\be
\label{eqn14}
\begin{array}{ll}
(A_1-\rho_1 I)^TX_1^{(k+1)}+X_1^{(k+1)}(A_1-\rho_1 I)-X_1^{(k+1)}S_1X_1^{(k+1)}\\
\\
\ds +\sum_{j=2}^N\delta_{1,j}Y_j^{(k)}+Q_1+2\rho_1Y_1^{(k)} \preceq 0.
\end{array}
\ee
In the meantime, we see by setting $i=1$ in (\ref{algo1}) that 
\be
\label{eqn15}
\begin{array}{ll}
(A_1-\rho_1 I)^TY_1^{(k+1)}+Y_1^{(k+1)}(A_1-\rho_1 I)-Y_1^{(k+1)}S_1Y_1^{(k+1)}\\
\\
\ds +\sum_{j=2}^N\delta_{1,j}Y_j^{(k)}+Q_1+2\rho_1Y_1^{(k)}=0.
\end{array}
\ee
Using Lemma \ref{lem2} and noting the stability of $A_1-\rho_1 I -S_1X_1^{(k+1)}$ from part (i) of Theorem \ref{thm1}, (\ref{eqn14}) and (\ref{eqn15}) lead to $X_1^{(k+1)} \succeq Y_1^{(k+1)}$.
\\

Next, suppose that for some $2 \le r \le N$, 
\be
\label{eqn16}
X_i^{(k+1)} \succeq Y_i^{(k+1)}, ~i=1, 2, \ldots, r-1.
\ee
It follows from (\ref{algo2}), (\ref{eqn13}), and (\ref{eqn16}) that 
\be
\label{eqn17}
\begin{array}{ll}
(A_r-\rho_r I)^TX_r^{(k+1)}+X_r^{(k+1)}(A_r-\rho_r I)-X_r^{(k+1)}S_rX_r^{(k+1)}\\
\\
\ds +\sum_{j=1}^{r-1}\delta_{r,j}Y_j^{(k+1)}+\sum_{j=r+1}^N\delta_{r,j}Y_j^{(k)}+Q_r+2\rho_rY_r^{(k)} \preceq 0.
\end{array}
\ee
Moreover, we see from (\ref{algo1}) and from the monotonicity of $\{Y_i^{(k)}\}$ established in Theorem \ref{thm3} that 
\be
\label{eqn18}
\begin{array}{ll}
(A_r-\rho_r I)^TY_r^{(k+1)}+Y_r^{(k+1)}(A_r-\rho_r I)-Y_r^{(k+1)}S_rY_r^{(k+1)}\\
\\
\ds +\sum_{j=1}^{r-1}\delta_{r,j}Y_j^{(k+1)}+\sum_{j=r+1}^N\delta_{r,j}Y_j^{(k)}+Q_r+2\rho_rY_r^{(k)} \succeq 0.
\end{array}
\ee
Using Lemma \ref{lem2} again and noting the stability of $A_r-\rho_r I -S_rX_r^{(k+1)}$ from part (i) of Theorem \ref{thm1}, (\ref{eqn17}) and (\ref{eqn18}) yield $X_r^{(k+1)} \succeq Y_r^{(k+1)}$, which implies that $X_i^{(k+1)} \succeq Y_i^{(k+1)}$, $i \in \langle N \rangle$.
\\

This finishes the proof.
\ep

Since with the assumptions of Theorems \ref{thm1} and \ref{thm3}, both Algorithms \ref{ricca} and \ref{accel} compute unique increasing sequences of positive semidefinite matrices, Theorem \ref{thm5} indicates that in this case, Algorithm \ref{accel} tends to converge faster than Algorithm \ref{ricca}.
\\

In the same spirit as Theorem \ref{thm5}, we can state below a parallel conclusion, whose proof is obvious and thus is omitted.

\begin{theorem}
\label{thm6}
Under the same assumptions as in Theorems \ref{thm2} and \ref{thm4} with $X_i^{(0)}=Y_i^{(0)}$ for all $i \in \langle N \rangle$, on letting $\{X_i^{(k)}\}$, $i \in \langle N \rangle$, be the sequences computed with Algorithm \ref{accel} and $\{Y_i^{(k)}\}$, $i \in \langle N \rangle$, be the corresponding sequences computed with Algorithm \ref{ricca}, we have that for each $i$, $X_i^{(k)} \preceq Y_i^{(k)}$, where $k=0, 1, 2, \ldots$.
\end{theorem}

The above Theorems \ref{thm5} and \ref{thm6}, together, provide an answer to the second open problem in \cite{Iva08}.
\\

Returning to the issue regarding the choice of $\rho_i$'s, similar to Algorithm \ref{ricca} --- see \cite[Remark 2]{CosVal04}, we now illustrate that these parameters should be picked in such a way that they are as small as possible. Numerical examples in this regard can be found in \cite{Iva08}.
\\

For the ease of statement, we first modify (\ref{algo2}) to that for $k=0, 1, 2, \ldots$, 
\be
\label{algo2+}
\begin{array}{ll}
[A_i-(\rho_i+\Delta \rho_i)I]^TY_i^{(k+1)}+Y_i^{(k+1)}[A_i-(\rho_i+\Delta \rho_i)I]-Y_i^{(k+1)}S_iY_i^{(k+1)}\\
\\
\ds +\sum_{j=1}^{i-1}\delta_{i,j}Y_j^{(k+1)}+\sum_{j=i+1}^N\delta_{i,j}Y_j^{(k)}+Q_i+2(\rho_i+\Delta \rho_i)Y_i^{(k)}=0,
\end{array}
\ee
where $i \in \langle N \rangle$ and $\Delta \rho_i \ge 0$ for all $i$; namely we consider a setting in which each $\rho_i$ in (\ref{algo2}) is augmented by $\Delta \rho_i$. Note that, here and in the sequel, we denote the sequences generated by (\ref{algo2+}) as $\{Y_i^{(k)}\}$ so as to differentiate them from $\{X_i^{(k)}\}$ generated by (\ref{algo2}). Clearly, the stabilizability and detectability conditions in Theorem \ref{thm1}, when it holds, still apply to (\ref{algo2+}).

\begin{theorem}
\label{thm7}
Under the same assumptions as Theorem \ref{thm1} with $X_i^{(0)}=Y_i^{(0)}$ for all $i \in \langle N \rangle$, let $\{X_i^{(k)}\}$, $i \in \langle N \rangle$, be the sequences computed from (\ref{algo2}) and let $\{Y_i^{(k)}\}$, $i \in \langle N \rangle$, be the corresponding sequences computed from (\ref{algo2+}), then we have that $X_i^{(1)} \succeq Y_i^{(1)}$, $i \in \langle N \rangle$.
\end{theorem}
\bp
Observe first that (\ref{algo2+}) satisfies all the conditions stated in Theorem \ref{thm1}. Hence, $\{Y_i^{(k)}\}$'s are uniquely determined by (\ref{algo2+}) and have all the properties in Theorem \ref{thm1}.
\\

Let us prove the conclusion by induction on $i$. At $i=1$, we obtain from $X_1^{(0)}=Y_1^{(0)} \preceq Y_1^{(1)}$ and (\ref{algo2+}) that 
\be
\label{eqn19}
\begin{array}{ll}
(A_1-\rho_1I)^TY_1^{(1)}+Y_1^{(1)}(A_1-\rho_1I)-Y_1^{(1)}S_1Y_1^{(1)}\\
\\
\ds +\sum_{j=2}^N\delta_{1,j}X_j^{(0)}+Q_1+2\rho_1X_1^{(0)}=2\Delta \rho_1(Y_1^{(1)}-Y_1^{(0)}) \succeq 0.
\end{array}
\ee
Comparing (\ref{eqn19}) and (\ref{algo2}) with $i=1$, and noting the stability of $A_1-\rho_1 I-S_1X_1^{(1)}$, we see $X_1^{(1)} \succeq Y_1^{(1)}$ by Lemma \ref{lem2}.
\\

Next, suppose that there exists some $2 \le r \le N$ such that $X_i^{(1)} \succeq Y_i^{(1)}$, $i=1, 2, \ldots, r-1$. This, together with $X_i^{(0)}=Y_i^{(0)} \preceq Y_i^{(1)}$ for all $i$ and (\ref{algo2+}), yield 
\be
\label{eqn20}
\begin{array}{ll}
\ds (A_r-\rho_rI)^TY_r^{(1)}+Y_r^{(1)}(A_r-\rho_rI)-Y_r^{(1)}S_rY_r^{(1)}+\sum_{i=1}^{r-1}\delta_{r,j}X_j^{(1)}\\
\\
\ds +\sum_{j=r+1}^N\delta_{r,j}X_j^{(0)}+Q_r+2\rho_rX_r^{(0)} \succeq 2\Delta \rho_r(Y_r^{(1)}-Y_r^{(0)}) \succeq 0.
\end{array}
\ee
Comparing (\ref{eqn20}) to (\ref{algo2}) with $i=r$, and in presence of the stability of $A_r-\rho_rI-S_rX_r^{(1)}$, we conclude using Lemma \ref{lem2} that $X_r^{(1)} \succeq Y_r^{(1)}$.
\\

Thus, $X_i^{(1)} \succeq Y_i^{(1)}$ for all $i \in \langle N \rangle$.
\ep

In Theorem \ref{thm7} above, for uniform satisfaction of the conditions in Theorem \ref{thm1} on both (\ref{algo2}) and (\ref{algo2+}), we follow \cite[Remark 2]{CosVal04} to perform only a ``single step'' analysis. This analysis, however, extends essentially to the scenario $X_i^{(k+1)} \succeq Y_i^{(k+1)}$, $i \in \langle N \rangle$, whenever $X_i^{(k)}=Y_i^{(k)}$ for all $i$. Accordingly, this result justifies that, in general, the larger $\rho_i$'s are, the slower the convergence (\ref{algo2}), i.e. Algorithm \ref{accel}, tends to exhibit.
\\

Finally, in the same vein as Theorem \ref{thm7}, we formulate here without proof its counterpart assuming the conditions in Theorem \ref{thm2}.

\begin{theorem}
\label{thm8}
Under the same assumptions as Theorem \ref{thm2} with $X_i^{(0)}=Y_i^{(0)}$ for all $i \in \langle N \rangle$, let $\{X_i^{(k)}\}$, $i \in \langle N \rangle$, be the sequences computed from (\ref{algo2}) and let $\{Y_i^{(k)}\}$, $i \in \langle N \rangle$, be the corresponding sequences computed from (\ref{algo2+}), then we have that $X_i^{(1)} \preceq Y_i^{(1)}$, $i \in \langle N \rangle$.
\end{theorem}

\section{Numerical Results}
\label{numer}

To illustrate our main conclusions in the preceding section, we present here relevant numerical results on one example. In accordance with the primary goals of this work, our numerical experiment has been carried out only with the Riccati iteration method, i.e. Algorithm \ref{ricca}, and the accelerated Riccati iteration method, i.e. Algorithm \ref{accel}. For numerical results comparing these methods with other existing methods, we refer the reader to \cite{Iva08}. Moreover, in view of our results, the example we provide here features distinct minimal and maximal positive semidefinite solutions.

\begin{example}
\label{examp}
Let $n=N=2$. Let 
$A_1=\left[\begin{array}{rr}
1 & -2\\
0 & -1
\end{array}\right]$, $A_2=\left[\begin{array}{cc}
1 & -1\\
0 & -3
\end{array}\right]$, $S_1=B_1B_1^T$, where $B_1=\left[\begin{array}{r}
5\\
-5
\end{array}\right]$, $S_2=B_2B_2^T$, where $B_2=\left[\begin{array}{r}
6\\
3
\end{array}\right]$, $\delta_{1,2}=2$, $\delta_{2,1}=3$, $Q_1=\left[\begin{array}{rr}
0 & 0\\
0 & 2
\end{array}\right]$, and $Q_2=\left[\begin{array}{cc}
0 & 0\\
0 & 3/2
\end{array}\right]$. Then, the CCARE as in (\ref{ccare}) has the minimal positive semidefinite solution $$X_1^-=\left[\begin{array}{cc}
  0.00000000 & 0.00000000\\
  0.00000000 & 0.28204532
\end{array}\right], ~X_2^-=\left[\begin{array}{cc}
  0.00000000 & 0.00000000\\
  0.00000000 & 0.27641488
\end{array}\right]$$
and the maximal positive semidefinite solution $$X_1^+=\left[\begin{array}{cc}
  0.50718185 & 0.24899225\\
  0.24899225 & 0.45594482
\end{array}\right], ~X_2^+=\left[\begin{array}{rr}
  0.32609148 & -0.16073063\\
 -0.16073063 & 0.48929635
\end{array}\right].$$
\end{example}

In our numerical experiment, the stopping criterion is set as $$\max_{i \in \langle 2 \rangle}\|X_i^{(k)}-X_i^{(k-1)}\|_F < {\rm tol}=10^{-8},$$ where $\| \cdot \|_F$ stands for the Frobenius norm. Upon the termination of either algorithm at the $m$-th iteration, the residual is calculated by $$\max_{i \in \langle 2 \rangle}\|{\cal R}_i(X_1^{(m)},X_2^{(m)})\|_F,$$ where ${\cal R}_i$ is given in (\ref{res}). In addition, for each $i$, we denote the largest eigenvalue of $X_i^{(k)}$ by $\lambda_1(X_i^{(k)})$, the smallest eigenvalue of $X_i^{(k)}$ by $\lambda_2(X_i^{(k)})$, and the spectrum of $X_i^{(k)}$ by $\sigma(X_i^{(k)})$, i.e. $\sigma(X_i^{(k)})=\{\lambda_1(X_i^{(k)}),\lambda_2(X_i^{(k)})\}$. These quantities are used in the illustrations.
\\

\begin{figure}[h!]
\begin{center}
\includegraphics[clip,scale=.48]{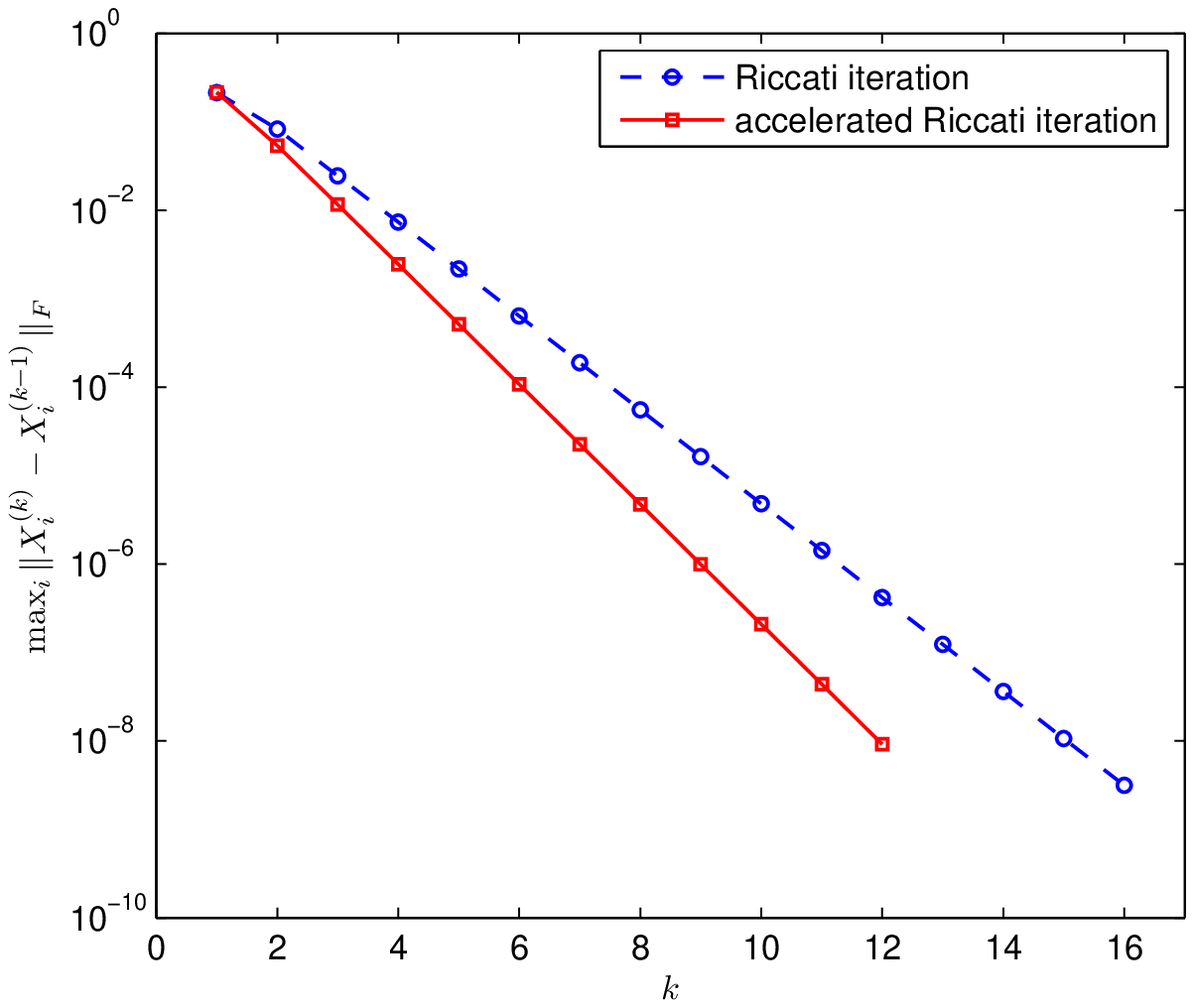} \hspace{-.3in} \includegraphics[clip,scale=.48]{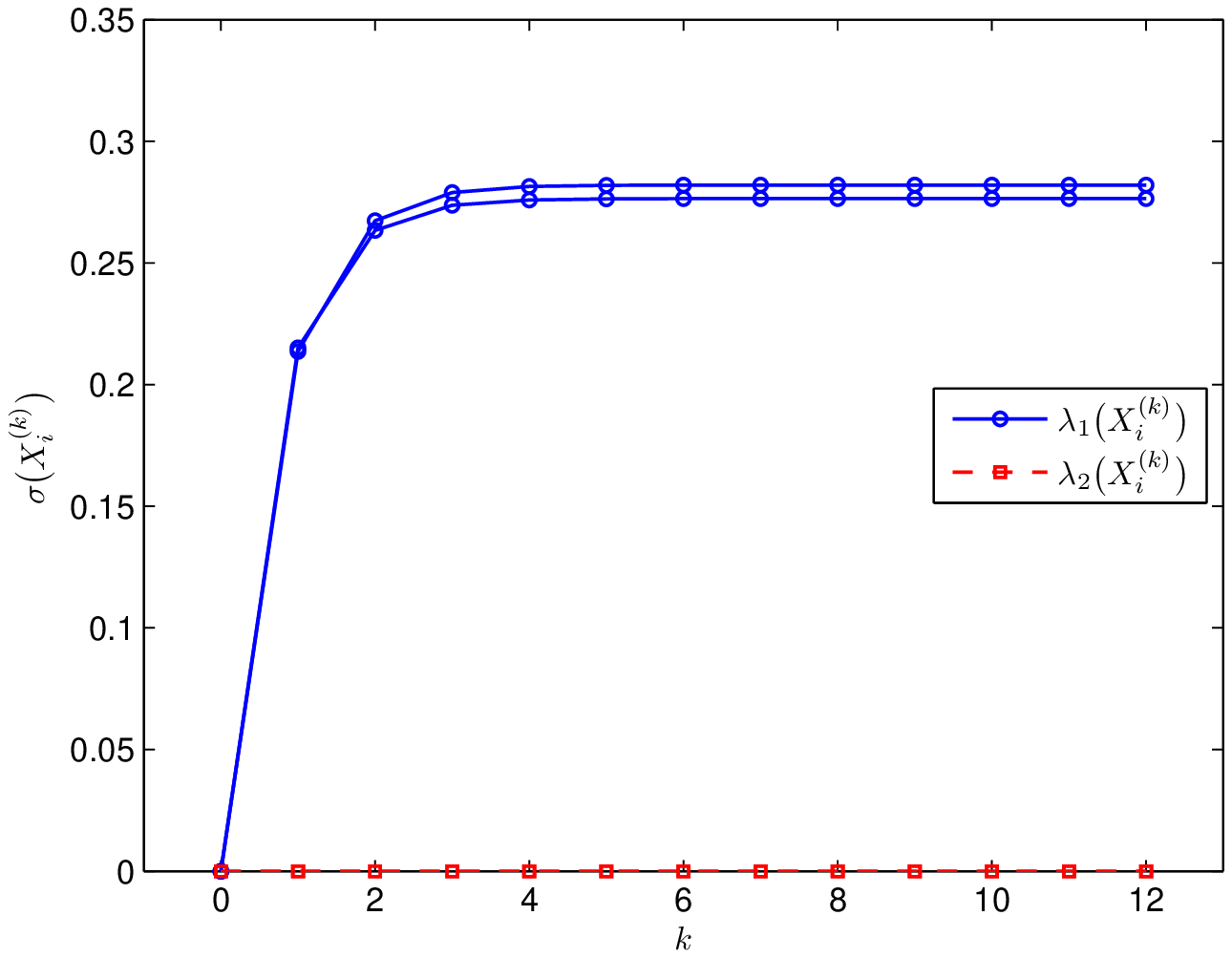}
\caption{The case of computing $X_i^-$ when $\rho_1=\rho_2=1.01$: The left panel shows $\max_{i \in \langle 2 \rangle}\|X_i^{(k)}-X_i^{(k-1)}\|_F$ from Algorithms \ref{ricca} and \ref{accel}, whereas the right panel illustrates the monotonic increasingness of the sequences $\{X_i^{(k)}\}$ obtained from Algorithm \ref{accel}. Note that in this case, $\lambda_2(X_i^{(k)})=0$ for all $i$ and all $k$.}
\label{fig1}
\end{center}
\end{figure}

To compute $X_i^-$, we choose $X_1^{(0)}=X_2^{(0)}=0$. It is not difficult to verify that the conditions of Corollary \ref{cor2} are all satisfied for any $\rho_i > 1$, $i=1, 2$. For $\rho_1=\rho_2=1.01$, Algorithm \ref{ricca} converges to $X_i^-$ in $16$ iterations, while Algorithm \ref{accel} does so in $12$ iterations, as shown in the left panel in Figure \ref{fig1}. In the meantime, the right panel of Figure \ref{fig1} displays the spectra of $X_i^{(k)}$ computed from Algorithm \ref{accel}, which shows that for each $i$, $\{X_i^{(k)}\}$ is monotonically increasing as confirmed by Corollary \ref{cor2}.
\\

\begin{table}[h!]
\begin{center}
\begin{tabular}{|c|c|c|c|c|}\hline
\multirow{2}{*}{$\rho_1=\rho_2$} & \multicolumn{2}{|c|}{Algorithm \ref{ricca}} & \multicolumn{2}{|c|}{Algorithm \ref{accel}}\\ \cline{2-5}
 & \# of iterations & residual & \# of iterations & residual \\ \hline
$1.5$ & $17$ & $3.92 \times 10^{-8}$ & $14$ & $4.25 \times 10^{-8}$\\
$1.1$ & $16$ & $1.91 \times 10^{-8}$ & $13$ & $1.43 \times 10^{-8}$\\
$1.01$ & $16$ & $1.20 \times 10^{-8}$ & $12$ & $3.48 \times 10^{-8}$\\ \hline
\end{tabular}
\caption{This table shows, as $\rho_i$ values vary, the numbers of iterations and residuals from Algorithms \ref{ricca} and \ref{accel} when computing $X_i^-$.}
\label{tab1}
\end{center}
\end{table}

With varying $\rho_i$ values, we summarize in Table \ref{tab1} the resulting numbers of iterations and residuals for computing $X_i^-$ by Algorithms \ref{ricca} and \ref{accel}. It points to that, as suggested by Theorem \ref{thm5}, Algorithm \ref{accel} converges faster than Algorithm \ref{ricca}. It also shows the speed-up in Algorithm \ref{accel} along with decreasing values of $\rho_i$, see Theorem \ref{thm7}.
\\

\begin{figure}[h!]
\begin{center}
\includegraphics[clip,scale=.48]{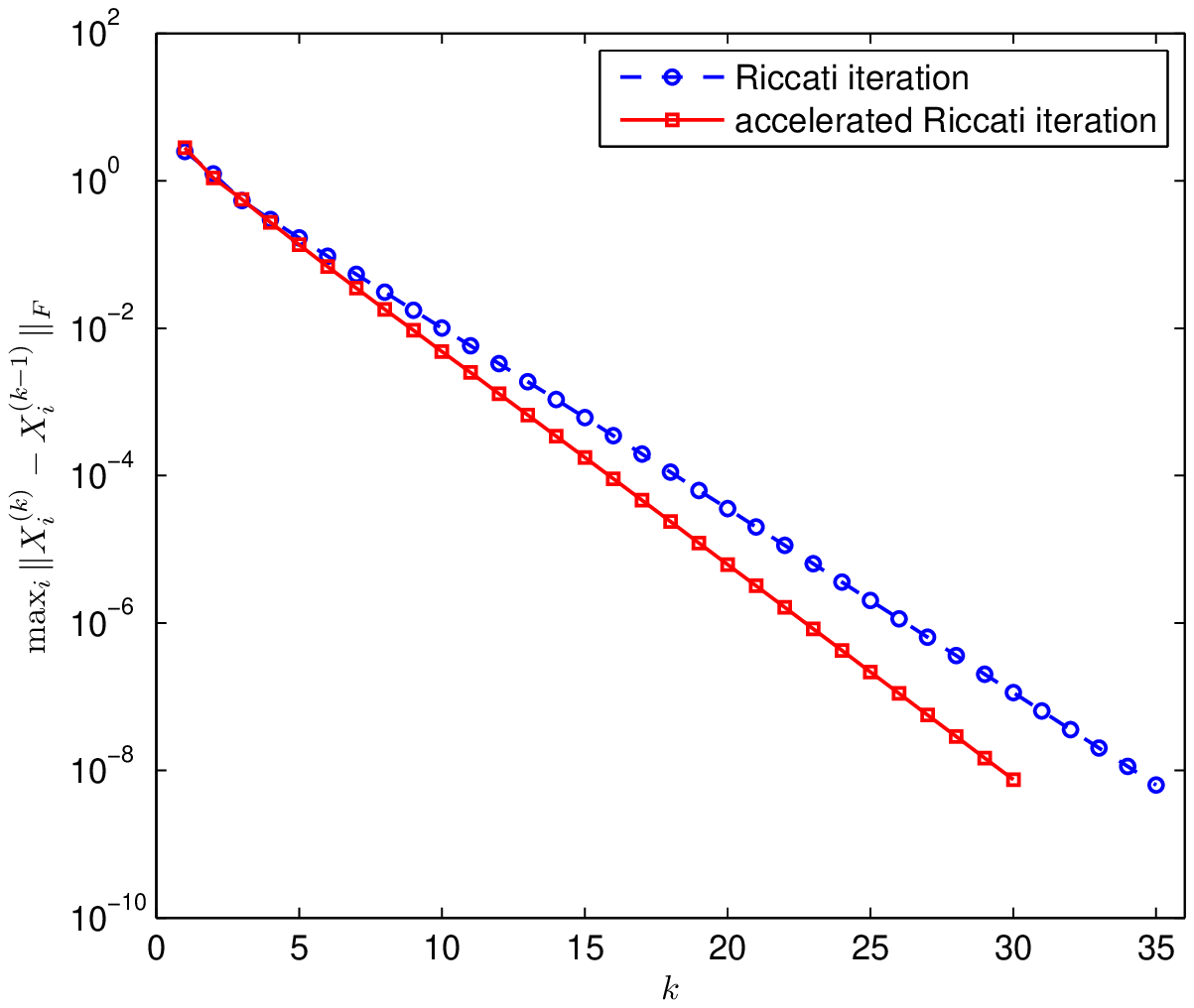} \hspace{-.3in} \includegraphics[clip,scale=.48]{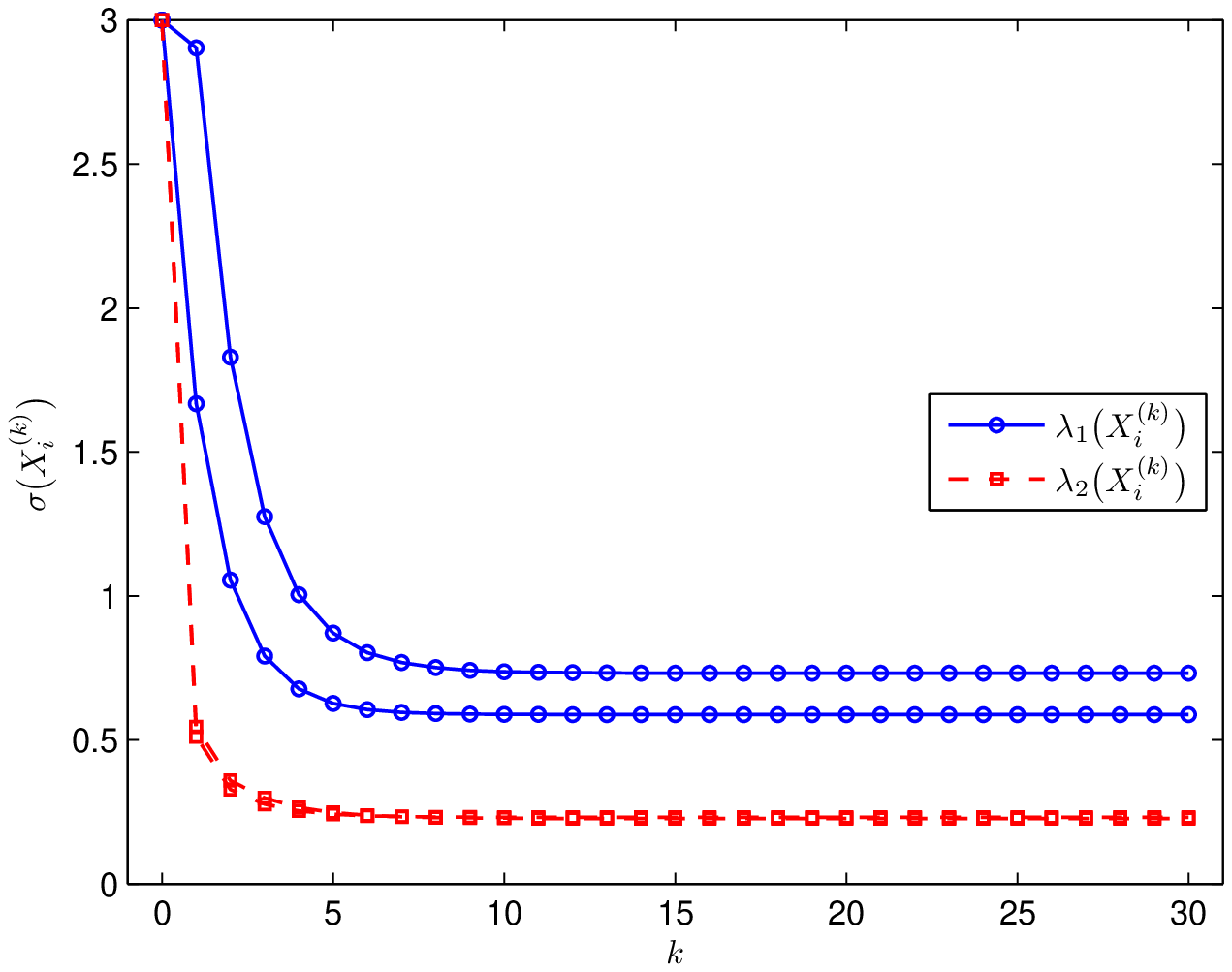}
\caption{The case of computing $X_i^+$ when $\rho_1=\rho_2=1.01$: The left panel shows $\max_{i \in \langle 2 \rangle}\|X_i^{(k)}-X_i^{(k-1)}\|_F$ from Algorithms \ref{ricca} and \ref{accel}, whereas the right panel shows the monotonic decreasingness of the sequences $\{X_i^{(k)}\}$ obtained from Algorithm \ref{accel}.}
\label{fig2}
\end{center}
\end{figure}

Next, to compute $X_i^+$, we choose $X_1^{(0)}=X_2^{(0)}=3I$. It is quite straightforward to verify that the conditions in Corollary \ref{cor4} are all satisfied for all $\rho_i > 1$, $i=1, 2$. Given $\rho_1=\rho_2=1.01$, Algorithm \ref{ricca} converges to $X_i^+$ in $35$ iterations, while Algorithm \ref{accel} does so in $30$ iterations, as illustrated by the left panel in Figure \ref{fig2}. In the right panel of Figure \ref{fig2}, the spectra of $X_i^{(k)}$ obtained from Algorithm \ref{accel} are plotted, showing that for each $i$, $\{X_i^{(k)}\}$ is monotonically decreasing.
\\

\begin{table}[h!]
\begin{center}
\begin{tabular}{|c|c|c|c|c|}\hline
\multirow{2}{*}{$\rho_1=\rho_2$} & \multicolumn{2}{|c|}{Algorithm \ref{ricca}} & \multicolumn{2}{|c|}{Algorithm \ref{accel}} \\ \cline{2-5}
 & \# of iterations & residual & \# of iterations & residual \\ \hline
$1.5$ & $42$ & $3.31 \times 10^{-8}$ & $38$ & $2.67 \times 10^{-8}$\\
$1.1$ & $36$ & $2.86 \times 10^{-8}$ & $32$ & $1.38 \times 10^{-8}$\\
$1.01$ & $35$ & $2.25 \times 10^{-8}$ & $30$ & $1.75 \times 10^{-8}$\\ \hline
\end{tabular}
\caption{This table shows, as $\rho_i$ values vary, the numbers of iterations and residuals from Algorithms \ref{ricca} and \ref{accel} in the case of computing $X_i^+$.}
\label{tab2}
\end{center}
\end{table}

With the same decreasing values of $\rho_i$ as in Table \ref{tab1}, we provide in Table \ref{tab2} evidence as indicated by Theorem \ref{thm8} of a speed-up in Algorithm \ref{accel} for computing $X_i^+$. As a comparison, the corresponding numerical results from Algorithm \ref{ricca} are given in Table \ref{tab2} as well. From these results, we also see that, as indicated by Theorem \ref{thm6}, Algorithm \ref{accel} tends to converge faster than Algorithm \ref{ricca} too when it comes to computing $X_i^+$.

\section{Concluding Remarks}
\label{concl}

The focus of this paper is on the two open problems raised in \cite{Iva08} concerning the monotone convergence of the accelerated Riccati iteration method as well as its rate of convergence in comparison with the pure Riccati iteration method. Our results aim mainly to settle these problems. In the process, we also broaden and strengthen some existing results in \cite{CosVal04}.
\\

A unique and quite useful feature of the Riccati iteration method and its accelerated version is their adoption of parameters $\rho_i$'s, which leads to easy satisfaction of the stabilizability and detectability conditions. In view of such parameters, we may call these methods ``shifted'' Riccati iteration methods as versus the ``unshifted'' Riccati iteration methods in \cite{ValGerCos99}.
\\

The idea of utilizing the updated $X_i^{(k+1)}$'s in the regular Riccati iteration method can be regarded as an extension to similar works on the accelerated Lyapunov iteration method \cite{Guo13, ValGerCos99}. These, besides \cite{Iva08}, have also motivated our development in this paper of theoretical results on the pure and accelerated Riccati iteration methods.
\\

Throughout this paper, we assume exact arithmetic in analyzing the two methods here. From a practical perspective, however, the stability and sensitivity analyses on these methods appear to be an interesting topic for future research.
\\

Another interesting topic for further investigation is a theoretical analysis comparing the performance of the two methods here with that of other existing numerical methods for solving the CCARE. In \cite{Iva08}, for example, we can find only numerical results concerning the performances of the methods under consideration here, Newton's method, together with the Lyapunov and the accelerated Lyapunov iteration methods. Nevertheless, several theoretical results on the performances of the ``unshifted'' Riccati methods and the Lyapunov iteration methods are presented in \cite{ValGerCos99}. We expect, therefore, that parallel results in this regard may also be developed to include the ``shifted'' Riccati iteration methods.
\\

Recalling the remark following Theorem \ref{thm4}, upper solution bounds play an important role in numerical computations on the CCARE. In fact, lower solution bounds are equally important. In Corollary \ref{cor2}, for example, $X_i^{(0)}$'s are indeed lower solution bounds. We feel that much work is still needed on simpler, tighter, and more easily applicable upper and lower solution bounds for the CCARE along with their applications in solving the CCARE numerically.
\\

Last but not least, the original framework in \cite{CosVal04} is more general in that it recasts the CCARE as one of the special cases from a so-called perturbed algebraic Riccati equation, abbreviated as PARE, involving a monotonically increasing positive semidefinite operator. It is one more important problem for us to explore as to whether the results here can be extended, with some splittings of that operator, to more effectively handle the general PARE.


\begin{thebibliography}{}

\bibitem{AboFreJan94}
H. Abou-Kandil, G. Freiling, G. Jank, Solution and asymptotic behavior of coupled Riccati equations in jump linear systems, IEEE Transactions on Automatic Control 39 (1994) 1631--1636.

\bibitem{BinLanMei12}
D. Bini, B. Lannazzo, B. Meini, Numerical Solution of Algebraic Riccati Equations, SIAM, Philadelphia, 2012. 

\bibitem{CosFraTod13}
O. Costa, M. Fragoso, M. Todorov, Continuous-Time Markov Jump Linear Systems, Springer-Verlag, Berlin, Heidelberg, 2013.

\bibitem{CosVal04}
E. Costa, J. do Val, An algorithm for solving a perturbed algebraic Riccati equation, European Journal of Control 10 (2004) 576--580.

\bibitem{CzoSwi01}
A. Czornik, A. Swierniak, Upper bounds on the solution of coupled algebraic Riccati equation, Journal of Inequalities and Applications 6 (2001) 373--385.

\bibitem{DamHin01}
T. Damm, D. Hinrichsen, Newton's method for a rational matrix equation occurring in stochastic control, Linear Algebra and Its Applications 332/334 (2001) 81--109.

\bibitem{DavShiWil08}
R. Davies, P. Shi, R. Wiltshire, Upper solution bounds of the continuous and discrete coupled algebraic Riccati equations, Automatica 44 (2008) 1088--1096.

\bibitem{GajBor95}
Z. Gajic, I. Borno, Lyapunov iterations for optimal control of jump linear systems at steady state, IEEE Transactions on Automatic Control 40 (1995) 1971--1975.

\bibitem{Guo13}
C. Guo, Iterative methods for a linearly perturbed algebraic matrix Riccati equation arising in stochastic control, Numerical Functional Analysis and Optimization 34 (2013) 516--529.

\bibitem{Iva07}
I. Ivanov, Iterations for solving a rational Riccati equation arising in stochastic control, Computers and Mathematics with Applications 53 (2007) 977--988.

\bibitem{Iva08}
I. Ivanov, On some iterations for optimal control of jump linear equations, Nonlinear Analysis 69 (2008) 4012--4024.

\bibitem{IvaHasMin01}
I. Ivanov, V. Hasanov, B. Minchev, On matrix equations $X\pm A^\ast X^{-2}A=I$, Linear Algebra and Its Applications 326 (2001) 27--44.

\bibitem{Mar90}
M. Mariton, Jump Linear Systems in Automatic Control, Marcel Dekker, New York, 1990.

\bibitem{ValGerCos99}
J. do Val, J. Geromel, O. Costa, Solutions for the linear-quadratic control problem of Markov jump linear systems, Journal of Optimization Theory and Applications 103 (1999) 283--311.

\bibitem{Wil71}
J. Willems, Least squares stationary optimal control and the algebraic Riccati equation, IEEE Transactions on Automatic Control 16 (1971) 621--634.

\bibitem{WilLaw07}
R. Williams II, D. Lawrence, Linear State-Space Control Systems, John Wiley \& Sons, Inc., Hoboken, New Jersey, 2007.

\bibitem{Xu13}
J. Xu, Unified, improved matrix upper bound on the solution of the continuous coupled algebraic Riccati equation, Journal of the Franklin Institute 350 (2013) 1634--1648.

\bibitem{XuRaj16}
J. Xu, P. Rajasingam, New unified matrix upper bound on the solution of the continuous coupled algebraic Riccati equation, Journal of the Franklin Institute 353 (2016) 1233--1247.

\bibitem{XuXia13}
J. Xu, M. Xiao, On the iterative refinement of matrix upper bounds for the solution of continuous coupled algebraic Riccati equations, Automatica 49 (2013) 2168--2175.

\end{thebibliography}
\end{document}